\documentclass[12pt,a4paper,leqno]{article}

\usepackage[top=1.1in, bottom=1.3in, left=1.5in, right=1.5in]{geometry}

\usepackage[ansinew]{inputenc}
\usepackage[T1]{fontenc}
\usepackage[english]{babel}
\usepackage{amsthm}
\usepackage{amsfonts}         
\usepackage{amsmath}
\usepackage{amssymb}
\usepackage{breqn}
\usepackage{bm}
\usepackage{url}

\newcommand{\R}{\mathbb{R}}
\newcommand{\C}{\mathbb{C}}

\newcommand{\N}{\mathbb{N}}

\newcommand{\Z}{\mathbb{Z}}

\theoremstyle{plain}
\newtheorem{theorem}{Theorem}
\newtheorem*{lemma*}{Lemma}
\newtheorem{lemma}{Lemma}

\newtheorem*{conj}{Conjecture}

\theoremstyle{remark}

\numberwithin{equation}{section}

\title{Applying Discrete Fourier Transform to the Hardy-Littlewood Conjecture}
\author{Jori Merikoski}
\date{\today}

\begin{document}

\maketitle

\begin{abstract}
We study the asymptotic behaviour of the prime pair counting function $\pi_{2k}(n)$ by the means of the discrete Fourier transform on $\Z/ n\Z$. The method we develop can be viewed as a discrete analog of the Hardy-Littlewood circle method. We discuss some advantages this has over the Fourier series  on $\R /\Z$, which is used in the circle method. We show how to recover the main term for $\pi_{2k}(n)$ predicted by the Hardy-Littlewood Conjecture from the discrete Fourier series. The arguments rely  on  interplay of Fourier transforms on $\Z/ n\Z$ and on its subgroup $\Z/ Q\Z,$  $Q \, \vert \, n.$

\end{abstract}

\tableofcontents
\section{Introduction}
In this article the main object of study is the counting function for prime number pairs
\begin{align}
\pi_{2k}(x)  := \vert \{ p \leq x: \, p \,\, \textrm{and} \,\, p+ 2k \,\, \textrm{both prime numbers} \} \vert,
\end{align}
where $2k \geq 2$ is an even integer. We use the notation $\vert A \vert$ for the cardinality of a given finite set $A.$ The methods used are very similar to the circle method, with the exception that we use Fourier analysis on $\Z/n\Z$ rather than $\R/\Z.$ Among studying $\pi_{2k}(x) ,$ the main purpose of this paper is to feature some advantages of using the discrete Fourier transform. The article is constructed as follows. In the first section we give a standard heuristic argument for the Hardy-Littlewood conjecture and prove Lemma \ref{easy}, which concerns a simplified situation. In the next section, we briefly recall some basic results on the discrete Fourier transform. In the third section, we use the Fourier transform  to study $\pi_{2k}(n).$ We manage to extract the main term predicted by the Hardy-Littlewood Conjecture from the Fourier series for $\pi_{2k}(n)$ by using Lemma \ref{easy}.   In the fourth section we give a brief analysis of the error terms. 

Let us begin with a well known heuristic argument for the asymptotic behaviour of $\pi_{2k}(x).$  Recall that if $\pi(x)$ denotes the number of prime numbers less than $x,$ then by the Prime Number Theorem $\pi(x) \sim x / \log x,$ as $ x \to \infty.$ That is, the density of prime numbers among the integers less than $x$ is $1/\log x .$ Thus, assuming that the event that some $n \leq x$ is a prime is independent of the number $n+2k$ being a prime, one might guess that $\pi_{2k}(x)$ is asymptotic to $x / \log^2 x.$ 

However, this assumption is flawed.  Let us first consider the case $2k=2.$ If $q>2$ is a prime we know already that $q+2$ is not divisible by 2. Thus $q+2$ is already twice more likely to be a prime than a random integer. Therefore, we should at least correct our guess $x / \log^2 x$ by multiplying by 2. Suppose then that $q$ is not divisible by a prime $p > 2.$ Then $q$ belongs to one of the $(p-1)$ non-zero residue classes of $p.$ If $q+2$ is also not divisible by $p,$ then it must be in one of the remaining $(p-2)$ of the $(p-1)$ non-zero residue classes of $p.$ Hence, given that $q$ is not divisible by $p,$ the number $q+2$ is not divisible by $p$ approximately $(p-2)/(p-1)$ of the time. Since the probability that the integer $q+2$ is not divisible by $p$ is $1-1/p,$ we should multiply our guess for $\pi_2(x)$ by a correction factor $\frac{p-2}{p-1} / \frac{p-1}{p} = \frac{p(p-2)}{(p-1)^2}.$ Doing this for all prime numbers $p > 2$ leads us to guess that
\begin{align}
\pi_2(x)\,  \sim \, 2 \left( \prod_{p> 2} \frac{p(p-2)}{(p-1)^2} \right)  \frac{x}{\log^2 x}.
\end{align} 
The convergence of the above product can be seen from
\begin{align}
\prod_{p> 2} \frac{p(p-2)}{(p-1)^2} =  \prod_{p> 2} \left( 1 - \frac{1}{(p-1)^2}\right).
\end{align}

 Let us then consider, for example, the cases $2k=2$ and $3k=6.$ If a number $q > 3$ is a prime, then it is not divisible by the number $3.$ This alone does not imply whether or not $q+2$ is divisible by $3,$ but we do know that $q+6$ also cannot be divisible by $3.$ This makes it more probable that $q+6$ is a prime than that $q+2$ is a prime. The number $q+2$ is not divisible by 3 only half of of the time, since either $q \equiv 1 \, (\textrm{mod} \, 3)$ or $q \equiv 2 \, (\textrm{mod} \, 3).$ Thus a reasonable guess is that $\pi_6(x) \sim 2 \, \pi_2(x).$

More generally, let $2k$ be a constant and let $q > k$ be a prime number.  Then  $q$ is not divisible by any of the prime factors of $2k.$ Hence also $q+2k$ is not divisible by any of the odd  prime factors of $k.$ Let $p$ be a odd prime factor of $k.$ Then $q+2$ is not divisible by $p$ approximately $1-1/(p-1)= \frac{p-2}{p-1}$ of the time, which leads us to guess that $\pi_{2k}(x) \sim   \prod_{2 < p \, \vert \, k} \frac{p-1}{p-2} \,\pi_2(x).$

This is part of the heuristic justification of the famous Hardy-Littlewood Conjecture on prime pairs (Hardy \& Littlewood, 1923). Hardy and Littlewood gave also other arguments to support the conjecture, using also the circle method from which the present article draws inspiration.

\begin{conj}\emph{\textbf{(Hardy-Littlewood Conjecture).}} \label{hl} Let $2k \geq 2$ be a constant. Then
\begin{align}
\pi_{2k}(x) \, \sim \, C_{2k} \frac{x}{\log^2 x}, 
\end{align}
as $x \to \infty,$ where the constant $C_{2k}$ is defined by
\begin{align}
C_{2k}: &= 2 \prod_{p > 2} \frac{p(p-2)}{(p-1)^2} \prod_{2 < p \, \vert \, k} \frac{p-1}{p-2} \\
& = \prod_{p \, \vert \, 2k } \frac{p}{p-1} \prod_{ p \,\nmid \, 2k } \frac{p(p-2)}{(p-1)^2} 
\end{align}
\end{conj}

Let us then look at a situation where the heuristics presented above can be made exact. Let $2k$ be an even positive integer and let $Q > 2k$ be any positive squarefree integer. As an easier version of the Hardy-Littlewood Conjecture, we can study how many of the integers $0 \leq a \leq Q-1$ satisfy $(a,Q) = 1 = (a+2k, Q),$ where $(a,b)$ is the greatest common divisor of $a$ and $b.$   If $Q$ is of the form $Q = \prod_{p<z} p,$ then the problem is an approximation of the Hardy-Littlewood Conjecture, since the prime numbers $> \sqrt{z}$ sit inside the set $\{a+kQ: \, (a,Q) = 1, \,  \, k \in \N  \}.$ We have the following lemma, proof of which is similar to the heuristic argument in the beginning. This result will be used later to study $\pi_{2k}(x)$.

\begin{lemma}\label{easy} If $Q$ is squarefree, then
\begin{align}
\sum_{\substack{ a= 0 \\ (a,Q) = 1 = (a+2k, Q)     }}^{Q-1} 1 = \prod_{p \, \vert \, (2k,Q)} (p-1) \prod_{\substack{p \, \vert \, Q \\ p \, \nmid  \,2k  }} (p-2)
\end{align}
\end{lemma}
\begin{proof}
The relations $(a,Q) = 1 = (a+2k, Q)$ hold precisely when for all prime numbers $p \, \vert \, Q,$ the prime $p$ does not divide $a$ or $a+2k.$ So let $p \, \vert \, Q$ be a prime, and let $d \, \vert \, Q$ be such that $p \, \nmid \, d$. If $p \, \vert \, 2k,$ then $p$ does not divide $a$ or $a+2k$ for  precisely $1-1/p$ of the numbers in $\{a: \, 0\leq a \leq Q-1,  (a,d) = 1 = (a+2k, d)\}$ If $p \, \nmid \, 2k,$ then $p$ does not divide $a$ or $a+2k$ for precisely $(1-1/p)(1-1/(p-1)) =   1-2/p$ of the numbers in $\{a: \, 0\leq a \leq Q-1, \,  (a,d) = 1 = (a+2k, d) \}.$  Multiplying over all primes $p \, \vert \, Q$ yields
\begin{align*}
\sum_{\substack{ a= 0 \\ (a,Q) = 1 = (a+2k, Q)     }}^{Q-1} 1 &= Q  \prod_{p \, \vert \, (2k,Q)} \left(1- \frac{1}{p}\right) \prod_{\substack{p \, \vert \, Q \\ p \, \nmid  \,2k  }} \left(1- \frac{2}{p}\right)  \\ 
 &=  \prod_{p \, \vert \, (2k,Q)} (p-1) \prod_{\substack{p \, \vert \, Q \\ p \, \nmid  \,2k  }} (p-2).
\end{align*}
\end{proof}

An alternative way to prove the above lemma is to use the inclusion-exclusion principle in a similar way as in the Legendre's Formula for $\pi(x);$ Let  $\nu_d(2k)$ be a completely multiplicative function with respect to the argument in the lower index $d,$ defined by setting $\nu_{p}(2k) = 1$ if $p \, \vert \, 2k,$ and $\nu_{p}(2k) = 2$ if $p \, \nmid\, 2k$  for primes $p.$ Then by the inclusion-exclusion principle
\begin{align*}
\sum_{\substack{ a= 0 \\ (a,Q) = 1 = (a+2k, Q)     }}^{Q-1} 1 &=  \sum_{d \, \vert \, Q} \mu(d) \frac{Q \, \nu_{d}(2k)}{d}  
 = Q  \prod_{p \, \vert \, Q} \left(1- \frac{\nu_p(2k)}{p}\right) \\
 &=  \prod_{p \, \vert \, (2k,Q)} (p-1) \prod_{\substack{p \, \vert \, Q \\ p \, \nmid  \,2k  }} (p-2),
\end{align*}
where $\mu$ is the M\"obius function.

\section{Fourier Transform Modulo $n$}
In this section we present some basic facts about the discrete Fourier transform. We use the notation $\Z_n := \Z/n\Z$ for the group of integers modulo $n.$ We identify $\Z_n$ with $\{ 1,2,3, \dots n\}$ as a subset of $\C$ in the formulas below. Let $f: \Z_n \to \C$ be a function. We may then define the discrete Fourier transform (modulo n) by
\begin{align}
\hat{f}(\xi) := \sum_{x \in \Z_n} f(x) e_n(-\xi x),
\end{align}
where $e_n( x) := e\left( \frac{x}{n}\right) = e^{\frac{2\pi i x}{n}}.$ 
We also use the notation $\mathcal{F}_n (f)(\xi)$ for the Fourier transform of $f$ modulo $n.$ We reserve the hat notation for Fourier transform  modulo $n.$ 

Many of the formulas of the usual Fourier transform hold also for the discrete one. We have collected some of them in the next lemma.
\begin{lemma} \emph{\textbf{(Fourier Transform In $\Z_n$).}} \label{four} For all functions $f,g: \Z_n \to \C $ we have
\begin{align}
f(x) &= \frac{1}{n} \sum_{\xi \in Z_n} \hat{f}(\xi) e_n(\xi x). \\
\sum_{\xi \in \Z_n} \hat{f}(\xi) \overline{\hat{g}(\xi)} &=  n \sum_{x\in \Z_n} 
 f(x)\overline{g(x)} \\
 \widehat{(f \ast g)}(\xi) &= \hat{f}(\xi)\hat{g}(\xi),
\end{align}
where the convolution $(f \ast g)$ is defined by
\begin{align*}
(f \ast g)(x) := \sum_{y \in \Z_n} f(y) g(x-y).
\end{align*}
\end{lemma}
\begin{proof}
The proofs of all three claims are direct computations based on a typical orthogonality relation . We have
\begin{align*}
\sum_{\xi \in Z_n} \hat{f}(\xi) e_n(\xi x) &= \sum_{\xi \in Z_n}  \sum_{y \in \Z_n}  f(y) e_n(\xi(x-y))  \\ 
&= \sum_{y \in \Z_n}  f(y) \sum_{\xi \in Z_n} e_n(\xi(x-y)) = n f(x),
\end{align*}
since
\begin{align*}
\sum_{\xi \in Z_n} e_n (\xi x) = \begin{cases} \frac{1-e_n( x n)}{1-e_n(x)} = 0, & x \neq 0, \\
 n, & x = 0. \end{cases}
\end{align*}
This proves the first formula.

For the discrete Plancherel's formula we have
\begin{align*}
\sum_{\xi \in \Z_n} \hat{f}(\xi) \overline{\hat{g}(\xi)} =   \sum_{x \in \Z_n} \sum_{y \in \Z_n}  f(x) \overline{g(y)} \sum_{\xi \in \Z_n}  e_n(\xi(y- x)) = n \sum_{x\in \Z_n} 
 f(x)\overline{g(x)} .
\end{align*}

To obtain the convolution formula we compute
\begin{align*}
 \widehat{(f \ast g)}(\xi) & = \sum_{x \in \Z_n} \sum_{y \in \Z_n} f(y) g(x-y) e_n(-\xi x)  \\
 &= \sum_{y \in \Z_n} f(y)e_n(-\xi y) \sum_{x \in \Z_n} g(x-y) e_n(-\xi (x-y)) = \hat{f}(\xi)\hat{g}(\xi).
\end{align*}
\end{proof}

We define the inverse Fourier transform
\begin{align*}
\mathcal{F}_n^{-1} (f) (x) := \frac{1}{n} \sum_{\xi \in Z_n} f (\xi) e_n(\xi x),
 \end{align*}
so that we have  $\mathcal{F}_n^{-1}(\mathcal{F}_n(f)) = f.$

\section{Applying the Discrete Fourier Transform to $\pi_{2k}(n)$}
Let $P:\N \to \{0,1\},$ denote the characteristic function for prime numbers. We use the same symbol $P$ for the function $\Z_n \to \{0,1\},$ which is given by the restriction of $P$ to $\{1,2,3, \dots n\}.$ For our purposes, the most useful property of the discrete Fourier transform is the convolution formula $ \widehat{(f \ast g)} = \hat{f}\hat{g},$ because we have
 \begin{align*}
 \pi_{2k}(n) = \sum_{x=1}^n P(x)P(x+2k) = (P \ast P_{-} ) (-2k) + \mathcal{O}(1),
 \end{align*}
where $P_{-}(x) := P(-x).$ We shall ignore the $\mathcal{O}(1)$ term, since we can redefine $ \pi_{2k}(n) := (P \ast P_{-} ) (-2k).$ Note that $\widehat{P_{-}} (\xi)= \overline{\hat{P}(\xi)}.$ Thus, by using the convolution formula and the inverse Fourier transform formula of Lemma \ref{four} we immediately obtain
\begin{theorem} \label{ptheur} Let $n$ and $k$ be positive integers. Then
\begin{align} \label{pheur}
\pi_{2k}(n) & = \frac{1}{n} \sum_{\xi \in \Z_n} \vert \hat{P}(\xi) \vert^2 \, e_n(-2k\xi).
\end{align}
In particular,
\begin{align}
\pi_{2k}(n) = \frac{n}{\log^2 n} +  \frac{1}{n} \sum_{\xi \in \Z_n \setminus \{ 0\}} \vert \hat{P}(\xi) \vert^2 \, e_n(-2k\xi) + o \left( \frac{n}{\log^2 n} \right), \quad n\to \infty.
\end{align}
\end{theorem} 
The second part  of the theorem follows from the Prime Number Theorem, since
\begin{align*}
\hat{P}(0) = \sum_{x \in \Z_n} P(x) = \pi(n) \sim \frac{n}{\log n}, \quad n \to \infty.
\end{align*} 
Theorem \ref{ptheur} is of great interest for the reason that the first term $\vert \hat{P}(0)\vert^2/n \sim n/\log^2 n$ is already of the right order of magnitude. We would like to show that that cancellations occur in the remaining  exponential sum. However, by the Hardy-Littlewood  Conjecture  we expect that the remaining sum is of the same order as the leading term. 

Taking just the first term in the sum \eqref{pheur} as an approximation to $\pi_{2k}(n)$ suffers from the same problem as the initial probably false heuristic guess $\pi_{2k}(n) \sim n/ \log^2 n$ that we gave in the beginning of this article; the primeness of an integer $p +2k$ is not independent of the number $p$ being a prime. That is, the set of prime numbers is expected to have some additive structure, which is described by the constants $C_{2k}$ in the Hardy-Littlewood Conjecture. For instance, if $p > 2$ is prime, it is odd. Then we already know that $p+2k$ is odd, and therefore twice as likely to be a prime number as a randomly chosen number. This argument gave the factor 2 in the constant $C_{2k}$. In the sum \eqref{pheur} this additive structure causes $\hat{P} (\xi)$ to be of order $n/\log^2 n$ for some of the $\xi \neq 0,$ as can be seen from the discussion below. This is an instance of a well known phenomenon that `additively random' sets should have small Fourier transforms for $\xi \neq 0$.

To see how to proceed from here, let us first look at how to recover the factor  2 of the Hardy-Littlewood Conjecture  from Theorem \ref{ptheur}. We may assume that $n$ is even, since if $n$ is odd we may consider $n+1.$ The error from this is clearly $\mathcal{O}(1).$ Then we have for all $0 \leq \xi < n/2 $
\begin{align*}
\hat{P}\left(\xi + \frac{n}{2}\right)  &= \sum_{x \in \Z_n} P(x) e_n \left(-\xi x - \frac{n}{2} x\right)   = \sum_{x \in \Z_n} P(x) e_n (-\xi x) e \left( - \frac{x}{2}\right) \\
& = -\hat{P}(\xi), 
\end{align*}
since $P(x) = 1$ only if $x$ is odd. Hence,
\begin{align*}
\pi_{2k}(n) & = \frac{1}{n} \sum_{\xi \in \Z_n} \vert \hat{P}(\xi) \vert^2 \, e_n(-2k\xi) = \frac{2}{n} \sum_{ 0 \, \leq \, \xi \, < \, \frac{n}{2} } \vert \hat{P}(\xi) \vert^2 \, e_n(-2k\xi) \\
& = 2 \frac{n}{\log^2 n} +  \frac{2}{n} \sum_{0 < \xi < \frac{n}{2} } \vert \hat{P}(\xi) \vert^2 \, e_n(-2k\xi) + o \left( \frac{n}{\log^2 n} \right).
\end{align*}
Note that the factor 2 was obtained using the fact that primes greater than 2 are odd, which corresponds to its heuristic justification.

To recover the whole constant $C_{2k},$ we require the Siegel-Walfisz Theorem on the distribution of primes in arithmetic progressions. The proof can be found in Davenport (1980).  It is similar to the proof of the Prime Number Theorem, except that it depends on the properties of $L$-functions instead of the $\zeta$-function. 
\begin{theorem}
Let 
\begin{align}
\pi(n,q,a) := \sum_{\substack{m \leq n \\ m \equiv a \mod q }} P(m),
\end{align}
and assume that the greatest common divisor $(a,q) = 1.$ Suppose that $q \leq \log^B n$ for some constant $B.$ Then there exists a constant $C_B$ depending only on $B$ such that
\begin{align}
\pi(n,q,a) = \frac{\text{\emph{Li}} (n)}{\varphi(q)} + \mathcal{O}\left( n \exp \left\{ - C_B \log^{\frac{1}{2}} n \right\} \right), \quad n \to \infty.
\end{align}
Here $\varphi(q)$ is the Euler tontient function, which gives the number of integers less than $q$ that are coprime to $q.$ $Li(n)$ is the offset logarithmic integral
\begin{align}
\text{\emph{Li}} (n) := \int_2^n \frac{dx}{\log x} \sim \frac{n}{\log n}, \quad n \to \infty .
\end{align} 
\end{theorem}

It should be noted that the next theorem does not give any estimate for the remaining sum. It only produces the factor $C_{2k}$ in a natural way from the exponential sum \eqref{pheur}. For the theorem we set
\begin{align*}
Q = Q_z := \prod_{p < z} p,
\end{align*}
where $z = z(n)$ is such that $Q_z \leq \log^B n.$ We also require that $z(n) \to \infty$ as $n \to \infty.$ We assume below that $n$ is divisible by $Q.$ If this is not the case, we can replace $n$ by $n+ a$ for some suitable $0 < a < Q.$ The error from this is clearly $\mathcal{O}(Q) = \mathcal{O}(\log^B n),$ so that we may restrict to the case $Q \, \vert \, n.$ In  the proof  we will use Fourier transform modulo $Q$ as well.
\begin{theorem} \label{coef} Let $k$ be a fixed positive integer, let $Q = Q_z$ be as above, and let $Q \, \vert \, n$. Then
\begin{align*}
\pi_{2k}(n) = C_{2k} \frac{n}{\log^2 n} + \frac{1}{n} \sum_{0 < \xi < \frac{n}{Q} }  e_n(-2k\xi) \sum_{r=0}^{Q-1} \left\vert \hat{P}\left(\xi +\frac{rn}{Q}\right) \right \vert^2  e \left( -\frac{2kr}{Q}\right) + o\left( \frac{n}{\log^2 n}\right),
\end{align*}
as $n \to \infty.$
\end{theorem}
\begin{proof}
Since $Q \, \vert \, n,$ we have by Theorem \ref{ptheur}
\begin{align*}
\pi_{2k}(n) &= \frac{1}{n}\sum_{\xi \in \Z_n} \vert \hat{P}(\xi) \vert^2 \, e_n(-2k\xi) = \frac{1}{n} \sum_{r=0}^{Q-1} \sum_{0 \leq \xi < \frac{n}{Q} }  \left\vert \hat{P}\left(\xi +\frac{rn}{Q}\right) \right \vert^2  e_n\left(-2k\xi -\frac{2krn}{Q}\right) \\
& = \frac{1}{n} \sum_{0 \leq \xi < \frac{n}{Q} } e_n\left(-2k\xi\right) \sum_{r=0}^{Q-1}  \left\vert \hat{P}\left(\xi +\frac{rn}{Q}\right) \right \vert^2  e \left( -\frac{2kr}{Q}\right) .
\end{align*}
To prove the theorem we need to show that for $\xi = 0$ we have
\begin{align*}
\frac{1}{n}\sum_{r=0}^{Q-1}  \left\vert \hat{P}\left(\frac{rn}{Q}\right) \right \vert^2  e \left( -\frac{2kr}{Q}\right)  \sim C_{2k}\frac{n}{\log^2 n}, \quad n \to \infty.
\end{align*}
We have
\begin{align}
\hat{P}\left(\frac{rn}{Q}\right) &= \sum_{x \in \Z_n} P(x)e \left( \frac{-xr}{Q}\right)  \nonumber \\ 
&=    \sum_{a= 0}^{Q-1}\left( \sum_{\substack{x \leq n \\ x \equiv a \mod Q } } P(x)\right) e\left( \frac{-ar}{Q}\right) \nonumber  \\
& = \sum_{a=0}^{Q-1} \pi(n,Q,a) e\left( \frac{-ar}{Q}\right). \label{here}
\end{align}
Hence, if we denote $\rho(a) := \pi(n,Q,a),$ then  $\hat{P}\left( rn/Q\right) = \mathcal{F}_Q(\rho) (r).$ Thus, by the Fourier inversion formula modulo $Q$ we have
\begin{align*}
\frac{1}{n}\sum_{r=0}^{Q-1}  \left\vert \hat{P}\left(\frac{rn}{Q}\right) \right \vert^2  e \left( -\frac{2kr}{Q}\right) & = \frac{Q}{n} \frac{1}{Q} \sum_{r=0}^{Q-1} \vert \mathcal{F}_Q(\rho) (r) \vert^2  e_Q(-2kr) \\
& = \frac{Q}{n} \mathcal{F}_Q^{-1} (F_Q(\rho \ast \rho_{-}))(-2k) \\
&= \frac{Q}{n} (\rho \ast \rho_{-})(-2k) \\
&= \frac{Q}{n} \sum_{r=0}^{Q-1} \pi(n,Q,r) \pi(n,Q,r+2k),
\end{align*}
where $\rho_{-}(x) := \rho(-x),$ and the convolution is modulo $Q.$ The contribution to the sum  from terms, where either $(r,Q) >1$ or $(r+2k,Q) >1,$ is clearly negligible. Thus, since the Siegel-Walfisz Theorem holds uniformly, the last expression is asymptotic to
\begin{align*}
\frac{\text{Li} (n) ^2}{n} \frac{Q}{\varphi(Q)^2} \left( \sum_{\substack{ a= 0 \\ (a,Q) = 1 = (a+2k, Q)     }}^{Q-1} 1 \right)  
\end{align*}
 since we can bound the summation of the error terms by
\begin{align*}
Q^3 \, n \exp \left\{ - C_B \log^{\frac{1}{2}} n \right\}  &\ll n \exp \left\{ 3 \log Q  - C_B \log^{\frac{1}{2}} n \right\}  \\
& \ll n \exp \left\{ 3 B \log \log n  - C_B \log^{\frac{1}{2}} n \right\} \\
& \ll n \exp \left\{ - C' \log^{\frac{1}{2}} n \right\}. 
\end{align*}
Hence, by Lemma \ref{easy} we have
\begin{align*}
 \frac{1}{n}\sum_{r=0}^{Q-1}  \left\vert \hat{P}\left(\frac{rn}{Q}\right) \right \vert^2  e \left( -\frac{2kr}{Q}\right) & \sim \frac{n}{\log^2 n} \frac{Q}{\varphi(Q)^2} \ \prod_{p \, \vert \, (2k,Q)} (p-1) \prod_{\substack{p \, \vert \, Q \\ p \, \nmid  \,2k  }} (p-2) \\ & = \frac{n}{\log^2 n} \prod_{p \vert Q} \frac{p}{(p-1)^2} \ \prod_{p \, \vert \, (2k,Q)} (p-1) \prod_{\substack{p \, \vert \, Q \\ p \, \nmid  \,2k  }} (p-2)  \\
& = \frac{n}{\log^2 n}  \prod_{p \, \vert \, (2k,Q) } \frac{p}{p-1} \prod_{\substack{ p \, \vert \, Q \\ p \,\nmid \, 2k }} \frac{p(p-2)}{(p-1)^2} \\
&\sim \frac{n}{\log^2 n}  \prod_{p \, \vert \, 2k } \frac{p}{p-1} \prod_{ p \,\nmid \, 2k } \frac{p(p-2)}{(p-1)^2}  = \frac{n}{\log^2 n}  C_{2k}, 
\end{align*}  
as $n \to \infty.$ This is because $Q = Q_z = \prod_{p < z} p$ tends to infinity as $n$ does and the second product is convergent. 
\end{proof}
The constant  $C_{2k}$ can be derived by a different method, which is similar to calculations that occur in the usual circle method, and which does not exploit the Fourier transform modulo $Q.$ We can use the Siegel-Walfisz Theorem already in (\ref{here}) to obtain
\begin{align*}
\hat{P}\left(\frac{rn}{Q}\right) = \frac{\text{{Li}} (n)}{\varphi(Q)} \sum_{\substack{1 \leq a < Q \\ (a,Q)= 1}} e\left( \frac{-ar}{Q}\right) + \mathcal{O}\left( Q \, n \exp \left\{ - C_B \log^{\frac{1}{2}} n \right\} \right).
\end{align*}
We now note that the sum above is just a Ramanujan's sum, often denoted by $c_Q(r)$. It can be evaluated as
\begin{align} \label{ram}
c_Q (r) = \mu \left( \frac{Q}{(r,Q)} \right) \frac{\varphi(Q)}{\varphi\left( \frac{Q}{(r,Q)}\right)}.
\end{align}
For proof of this, see Davenport (1980, p. 148). This leads us to the expression
\begin{align*}
\frac{1}{n}\sum_{r=1}^{Q}  \left\vert \hat{P}\left(\frac{rn}{Q}\right) \right \vert^2 e \left( -\frac{2kr}{Q}\right)  \sim \frac{n}{\log^2 n} \sum_{r=1}^{Q}    \left \vert \mu \left( \frac{Q}{(r,Q)} \right) \right\vert \frac{1}{\varphi\left( \frac{Q}{(r,Q)}\right)^2} e \left( -\frac{2kr}{Q}\right) ,
\end{align*}
The sum can be evaluated by writing the sum as running over the divisors of $Q.$ This can be done since  $Q/(r,Q)$ always divides $Q.$ That is, any $d \, \vert \, Q$ appears in the  sum whenever $Q/(r,Q) = d.$ This equivalent to  that $r$ is of the form $r = Qb/d,$ where $(b,d) = 1$ and $b \leq d.$ Therefore, collecting all the terms where $Q/(r,Q) = d$ together yields
\begin{align*}
\sum_{r=1}^{Q}    \left \vert \mu \left( \frac{Q}{(r,Q)} \right) \right\vert \frac{1}{\varphi\left( \frac{Q}{(r,Q)}\right)^2} e \left( -\frac{2kr}{Q}\right) &= \sum_{d \, \vert \, Q} \frac{\vert \mu(d) \vert}{\varphi(d)^2} \sum_{\substack{b \, \leq \, d (b,d) = 1}} e \left( -\frac{2kb}{d}\right) \\
&= \sum_{d \, \vert \, Q} \frac{\vert \mu(d) \vert}{\varphi(d)^2} c_d(2k).
\end{align*}
The last expression is reminiscent of the so called `singular series,' which were analyzed all the way back in  Hardy \& Littlewood (1923). 
The last sum is equal to
\begin{align*}
\sum_{d \, \vert \, Q} \vert \mu(d) \vert \prod_{p \, \vert \, d} \frac{1}{(p-1)^2}  c_d( 2k)
& = \prod_{p \, \vert \, Q} \left(1+\frac{c_p(2k)}{(p-1)^2} \right) \\
& = \prod_{p \, \vert \, (2k,Q) } \frac{p}{p-1} \prod_{\substack{ p \, \vert \, Q \\ p \,\nmid \, 2k }} \frac{p(p-2)}{(p-1)^2} \\
&\to \prod_{p \, \vert \, 2k } \frac{p}{p-1} \prod_{ p \,\nmid \, 2k } \frac{p(p-2)}{(p-1)^2}  = C_{2k}, 
\end{align*}  
as $n \to \infty.$ This is because $Q = Q_z = \prod_{p < z} p$ tends to infinity as $n$ does and the second product is convergent.  In the above we have used the  multiplicativity of $d \mapsto c_d(2k)$ to obtain the first equality (Davenport, 1980, p. 149). The second equality holds, because by \eqref{ram} we have
\begin{align*}
c_p (2k) = \begin{cases} p-1, & p \, \vert \, 2k, \\
-1, & p \nmid 2k.
\end{cases}
\end{align*}
This gives a second proof for the claim. By using Fourier transform modulo $Q$ in the first proof, we were able to avoid  the singular series and Ramanujan's sums entirely.

As a quick remark, we note that all the calculations done in this section for $\pi_{2k}(n)$ apply also to the function
\begin{align*}
\psi_{2k}(n) := \sum_{x \in \Z_n} \Lambda(x)\Lambda(x+2k) = (\Lambda \ast \Lambda_{-})(-2k).
\end{align*}
It is sometimes more convenient to use this function because of the nice properties of the von Mangoldt function. It is a straightforward computation to see that the Hardy-Littlewood Conjecture is equivalent to $\psi_{2k}(n) \sim C_{2k} n.$ The Siegel-Walfisz Theorem holds also for the function
\begin{align}
\psi(n,q,a) := \sum_{\substack{m \leq n \\ m \equiv a \mod q }} \Lambda(m),
\end{align}
in the form
\begin{align}
\psi(n,q,a) = \frac{x}{\varphi(q)} + \mathcal{O}\left( n \exp \left\{ - C_B \log^{\frac{1}{2}} n \right\} \right), \quad n \to \infty
\end{align}
for $q \leq \log^B n$ and  $(a,q)= 1$ (Davenport, 1980). Therefore, the same calculations as for the function $\pi_{2k}(n)$ give us the following theorem.
\begin{theorem} Let $k$ be fixed, and let  $Q$ be as in Theorem \ref{coef}. Then
\begin{align}
\psi_{2k}(n) & = \frac{1}{n} \sum_{\xi \in \Z_n} \vert \hat{\Lambda}(\xi) \vert^2 \, e_n(-2k\xi),
\end{align}
and
\begin{align}
\psi_{2k}(n) = C_{2k} n + \frac{1}{n} \sum_{0 < \xi < \frac{n}{Q} }  e_n(-2k\xi) \sum_{r=0}^{Q-1} \left\vert \hat{\Lambda}\left(\xi +\frac{rn}{Q}\right) \right \vert^2  e \left( -\frac{2kr}{Q}\right) + o\left( n \right),
\end{align}
as $n \to \infty.$
\end{theorem}

\section{Analysis of the Error Term}
In this section we look into the more difficult problem of estimating the remaining terms. We  note that immediately we have some control over the sizes of $\hat{P}(\xi).$ The trivial estimate is of course
\begin{align*}
\hat{P}(\xi) = \sum_{x \in \Z_n} P(x) e_n(-\xi x) = \mathcal{O}\left( \frac{n}{\log n}\right).
\end{align*} 
The discrete analog of the Plancherel's formula (second part of Theorem \ref{four}) gives us much stronger constraints. In particular, the next theorem implies that $\frac{1}{n} \vert \hat{P}(\xi)\vert^2$ can be of size $n/C \log^2 n$ for at most $C \log n$ of $\xi \in \Z_n.$
\begin{theorem} Let $n$ be a positive integer. Then
\begin{align}
\frac{1}{n} \sum_{\xi \in \Z_n} \vert \hat{P}(\xi) \vert^2 = \sum_{x \in \Z_n} P(x) \sim \frac{n}{\log n}, \quad n \to \infty.
\end{align}
\end{theorem}

Even so, we expect that to bound the error term, we would have to take into account delicate cancellations coming from the exponentials in Theorem \ref{coef}. Motivated by the form of the theorem, and by our success in using the Fourier transform modulo $Q$ for the main term, we suggest that the error term 
\begin{align*}
\frac{1}{n} \sum_{0 < \xi < \frac{n}{Q} }  e_n(-2k\xi) \sum_{r=0}^{Q-1} \left\vert \hat{P}\left(\xi +\frac{rn}{Q}\right) \right \vert^2  e \left( -\frac{2kr}{Q}\right)
\end{align*}
should first be bounded by the triangle inequality as
\begin{align*}
\frac{1}{Q} \max_{\xi = 1,2,\dots, \frac{n}{Q}}\left \vert  \sum_{r=0}^{Q-1} \left\vert \hat{P}\left(\xi +\frac{rn}{Q}\right) \right \vert^2  e_Q(-2kr)  \right \vert= \max_{\xi = 1,2,\dots, \frac{n}{Q}} \left \vert  \mathcal{F}_Q^{-1} (\vert \mathcal{F}_n(P)_{\xi} \vert^2) (-2k) \right \vert,
\end{align*}
where $\mathcal{F}_n(P)_{\xi}$ denotes the function $\Z_Q \to \C,$ $a \mapsto \hat{P}\left(\xi +an/Q\right).$ This form emphasizes the interaction between  Fourier transforms on $\Z_n$ and on its subgroup $\Z_Q$. Following the proof of Theorem \ref{coef}, we obtain
\begin{align*}
\hat{P}\left(\xi + \frac{rn}{Q}\right) &= \sum_{x \in \Z_n} P(x)e_n(-x\xi) e \left( \frac{-xr}{Q}\right)   \\ 
&=    \sum_{a= 0}^{Q-1}\left( \sum_{\substack{x \leq n \\ x \equiv a \mod Q } } P(x) e_n(-x\xi) \right) e_Q(-xr) \\
& = :\sum_{a=0}^{Q-1} \pi_{\xi}(n,Q,a) e_Q(-xr)
\end{align*}
Hence, if we denote $\rho_{\xi}(a) := \pi_{\xi}(n,Q,a),$ then  $\hat{P}\left( \xi + rn/Q\right) = \mathcal{F}_Q(\rho_{\xi}) (r).$ Thus, by the Fourier inversion formula modulo $Q$ we have
\begin{align*}
\mathcal{F}_Q^{-1} (\vert \mathcal{F}_n(P)_{\xi} \vert^2) (-2k) 
& =  \mathcal{F}_Q^{-1} (F_Q(\rho_{\xi} \ast (\rho_{\xi})_{-}))(-2k) \\
&= (\rho_{\xi} \ast (\rho_{\xi})_{-}))(-2k) \\
&= \sum_{a=0}^{Q-1} \pi_{\xi}(n,Q,a) \pi_{\xi}(n,Q,a+2k),
\end{align*}
where
\begin{align}
& \pi_{\xi}(n,Q,a) = \left( \sum_{\substack{x \leq n \\ x \equiv a \mod Q } } P(x) e_n(-x\xi) \right)  \label{estimate1}\\ 
& \pi_{\xi}(n,Q,a+2k) = \left( \sum_{\substack{-2k < x \leq n -2k \\ x \equiv a \mod Q } } P(x+2k) e_n(-(x+2k)\xi) \right). \label{estimate2}
\end{align}
Thus the error term is bounded by
\begin{align*}
\max_{\xi = 1,2,\dots, \frac{n}{Q}} \left \vert \sum_{a=0}^{Q-1} \pi_{\xi}(n,Q,a) \pi_{\xi}(n,Q,a+2k)\right \vert
\end{align*}
In effect, we have reduced bounding the error term to showing that $\pi_{\xi}(n,Q,a)$ and $\pi_{\xi}(n,Q,a+2k)$ cannot simultaneously be very large. By the Siegel-Walfisz Theorem, the maximal possible order for $\pi_{\xi}(n,Q,a)$ is $n/(\varphi(Q)\log n).$ Since $Q =\prod_{p<z} p,$ we have $\varphi(Q) = Q \prod_{p<z} (1-1/p) \sim e^{-\gamma} Q /\log z$ by Merten's Theorem (Hardy \& Wright, 1945, p. 349). Hence, for $Q \sim \log^B x$  we have a trivial estimate $\pi_{\xi}(n,Q,a) = \mathcal{O} (n \log z/\log^{B+1} n).$ Since $\pi_{\xi}(n,Q,a)$ is defined as a Fourier transform modulo $n$ we have by the discrete Plancherel's formula
\begin{align*}
\frac{1}{n} \sum_{\xi \in \Z_n} \vert \pi_{\xi}(n,Q,a) \vert^2 = \sum_{\substack{x \leq n \\ x \equiv a \mod Q } } P(x) \sim \frac{1}{\phi(Q)} \frac{n}{\log n}, \quad n \to \infty.
\end{align*}
Hence, the typical $\pi_{\xi}(n,Q,a)$ is of order $\sqrt{n/ (\varphi(Q)\log n)},$ so that the typical product  $\pi_{\xi}(n,Q,a)\pi_{\xi}(n,Q,a+2k)$ is of order $n/ (\varphi(Q)\log n) = \mathcal{O} (n \log z/\log^{B+1} n),$ which just  logarithmic factors larger than what we want. Thus, it might be possible to get the neccessary savings by showing that $\pi_{\xi}(n,Q,a)$ and $\pi_{\xi}(n,Q,a+2k)$ cannot simultaneously be large.

\section{Conclusions}
There are many benefits from using the discrete Fourier transform instead of Fourier series on $\R/\Z$ which is commonly used. Firstly, one avoids the technicalities involved in defining the major and minor arcs. Secondly, by utilizing Fourier transform modulo $Q,$ we easily avoid having to deal with singular series, which can be cumbersome, as can be seen from from the alternative proof of Theorem \ref{coef}. The argument becomes more rigid in nature, and it emphasizes the algebraic ideas underlying the computations. Finally, extracting the main term from the Fourier series follows the heuristic argument behind the Hardy-Littlewood Conjecture, making calculations more intuitive. 

We have shown that
\begin{align*}
\pi_{2k}(n) &= C_{2k} \frac{n}{\log^2 n} + \frac{1}{n} \sum_{0 < \xi < \frac{n}{Q} }  e_n(-2k\xi) \sum_{r=0}^{Q-1} \left\vert \hat{P}\left(\xi +\frac{rn}{Q}\right) \right \vert^2  e \left( -\frac{2kr}{Q}\right) + o\left( \frac{n}{\log^2 n}\right) \\
&=  C_{2k} \frac{n}{\log^2 n}  + \mathcal{O}\left( \max_{\xi = 1,2,\dots, \frac{n}{Q}} \left \vert    \sum_{a=0}^{Q-1} \pi_{\xi}(n,Q,a) \pi_{\xi}(n,Q,a+2k) \right \vert\right)  + o\left( \frac{n}{\log^2 n}\right).
\end{align*}

Possible topic for future papers would be to  apply the method to related problems, for example,  to give a proof for Vinogradov's Three Primes Theorem using the discrete Fourier transform in a similar fashion as above.

\section*{References}
\addcontentsline{toc}{section}{References}
Davenport, H. (1980). \emph{Multiplicative Number Theory.} Second edition. Springer, New York. \\ \\
Hardy, G.  \&  Littlewood,  J.  (1923). Some Problems of 'Partitio Numerorum.' III. On the Expression of a Number as a Sum of Primes. \emph{Acta Mathematica,} 44, pp. 1-70. 
\\ \\
Hardy, G. \& Wright, E. (1945). \emph{An Introduction to the Theory of Numbers.} Second edition. Oxford University Press, Oxford.
\end{document}